\documentclass[12pt]{book}

\usepackage{amsthm}
\usepackage{amsmath}
\usepackage{amscd}
\usepackage{amssymb}

\usepackage{graphicx}
\DeclareGraphicsExtensions{mps,pdf}
\usepackage{wrapfig}

\ifx\pdfoutput\unknown\fi

\begin{document}

\newcommand{\pd}{\partial}
\newcommand{\R}{\mathbb R}
\newcommand{\C}{\mathbb C}
\newcommand{\Sb}{\mathbb S}
\newcommand{\D}{\mathbb D}
\newcommand{\Z}{\mathbb Z}
\newcommand{\N}{\mathbb N}

\newcommand{\spt}{\operatorname{supp}}
\newcommand{\Hom}{\operatorname{Hom}}
\newcommand{\Ker}{\operatorname{Ker}}
\newcommand{\id}{\operatorname{id}}
\newcommand{\Der}{\operatorname{Der}}
\newcommand{\sh}{\operatorname{sh}}
\newcommand{\ch}{\operatorname{ch}}
\newcommand{\tnh}{\operatorname{tnh}}
\newcommand{\ctnh}{\operatorname{ctnh}}
\newcommand{\ctg}{\operatorname{ctg}}
\newcommand{\tg}{\operatorname{tg}}
\newcommand{\tgh}{\operatorname{tgh}}

\centerline{{\sc Mechanical Normal Forms of Knots and Flat Knots}}

\medskip  
\centerline{{\sc A.B.Sossinsky}} 

\medskip
\centerline{Independent University of Moscow}

\centerline{Institute for Problems in Mechanics RAS} 
\bigskip

\begin{flushright}
{\it Dedicated to the memory of V.I.Arnold}
\end{flushright}

\small
{\bf Abstract.} A new type of knot energy is presented via real life experiments involving a thin resilient
metallic tube. Knotted in different ways, the device mechanically acquires a uniquely determined (up to
isometry) normal form at least when the original knot diagram has a small number of crossings, thus
outperforming the famous M\"obius energy due to Jun O'Hara and studied by Michael Freedman et al.
Various properties of the device are described (under certain conditions it does the Reidemeister and
Markov moves, it beautifully performs the Whitney trick by uniformizing its own local curvature). If the
device is constrained between two parallel planes (e.g. glass panes), it yields a real life model of a flat knot
(class of knot diagrams equivalent under Reidemeister $\Omega_2$ and $\Omega_3$ moves) also leading to
uniquely determined "flat normal forms" (for a small number of crossing points of the given flat knot diagram).
The paper concludes with two mathematical theorems, one reducing the knot recognition problem to the
flat knot recognition problem, the other (due to S.V.Matveev) giving an easily computable complete system
of invariants for the flat unknot knot equivalence problem. 
\normalsize

\bigskip
In one of his lectures, Vladimir Arnold declared that ``mathematics is that part of physics  where experiments are cheap". The present paper is an illustration of this thesis.  

Accordingly, the paper begins (Sec.1) with a description of our ``cheap experiments", which are performed with a thin flexible but resilient wire. This wire can be knotted and placed on a horizontal table, and in all the experiments  it moves in a very specific way, sometimes sliding along the table top, sometimes jumping up in the air, and reaches an equilibrium position, which we call the mechanical normal form of the knot. Surprisingly, the experiments show that, for a small number of crossings of the initial knot diagram, the normal form is always unique up to isometry. In this situation, it turns out that our simple device performs much better than the celebrated M\"obius energy models [1], [2] in computer simulations (see, e.g. [3] or [4]).

In Sec.2, we discuss possible idealizations of our physical model of fat knots and define 
what we call the $s$-normal form of a classical knot, where $s>0$ is a very small parameter and the knot is understood as a smooth curve 
rather than a ``fat knot" or a tube. 

In Sec.3, we describe the situations in which our device jumps off the table in terms of the Whitney index (winding number) of plane curves.
 
In Sec.4, we analyze the passage to the normal form and note that our clever device can perform all the Reidemeister moves [5, p.11], the Markov moves [5, p.60], and the Whitney trick [6, p.17] in a very impressive way. 

In Sec.5, we discuss the nonuniqueness of normal forms and present conjectural candidates for counterexamples to uniqueness. 

In Sec.6, we consider flat knots (i.e. knots lying on the plane except in the vicinity of crossing points, where they rise a bit above the plane, eliminating the self-intersections 
of the vertical projection). Flat knots may be regarded as purely mathematical objects,  we discuss three of their invariants and reduce the verification of the equivalence of knot diagrams (as curves in space) to the verification of ``flat equivalence" of flat knots (plane curves with underpass-overpass information).

In the concluding Sec.7, we (very optimistically) analyze the perspectives of this research, asserting that an energy functional may be devised to mimic  our mechanical device, and this will allow us to continue our experiments electronically (still very cheap!). We conjecture that our approach will yield a fast algorithm for comparing knots with a small number of crossings presented by their diagrams. Further, we claim that our notion of  $s$-normal form allows to define new knot energies (not necessarily mimicking our device) and involving, perhaps, the twisting number, and/or the torsion of the smooth curve defining the knot, as well as other characteristics of the knot.

\bigskip
{\bf 1. The mechanical model}

\smallskip
Our mechanical model of fat knots is constructed from a flexible but resilient  hollow rectilinear metallic cylinder of length 50 cm and diameter 2mm. If one bends the cylinder and places it on a table, it straightens out almost instantly (after some oscilliations) and reacquires its original rectilinear shape; if one bends it so that several loops are formed and the cylinder is released, it will jump up slightly off the table and  straighten out as it falls back on the table again, having reacquired its rectilinear form. We can say that, when released, it minimizes its total curvature as well as its curvature at each point, almost instantly making the curvature equal to zero everywhere.   

One extremity of the cylinder ends in a spike which fits snugly into the hollow at the other  
extremity, so that one can bend it, forming a plane closed Jordan curve; if released, it then forms a perfectly round circle lying on the table. One can also bend it into a nontrivial knot diagram; placed on the table, the mechanical knot immediately changes shape (sometimes sliding along the table, sometimes first jumping vertically up) and stabilizes in a position that we call its {\it mechanical normal form}. In this form, our device lies flat on the table, only rising above it by 4 mm at the overpasses. We say that the knot is {\it almost planar} or {\it flat}.

How unique is this normal form?    
Having performed thousands of experiments in which we tied the cylinder into different knots (with a small number of crossings), we obtained the following amazing result. 

\begin{figure}[h]\centering
\includegraphics{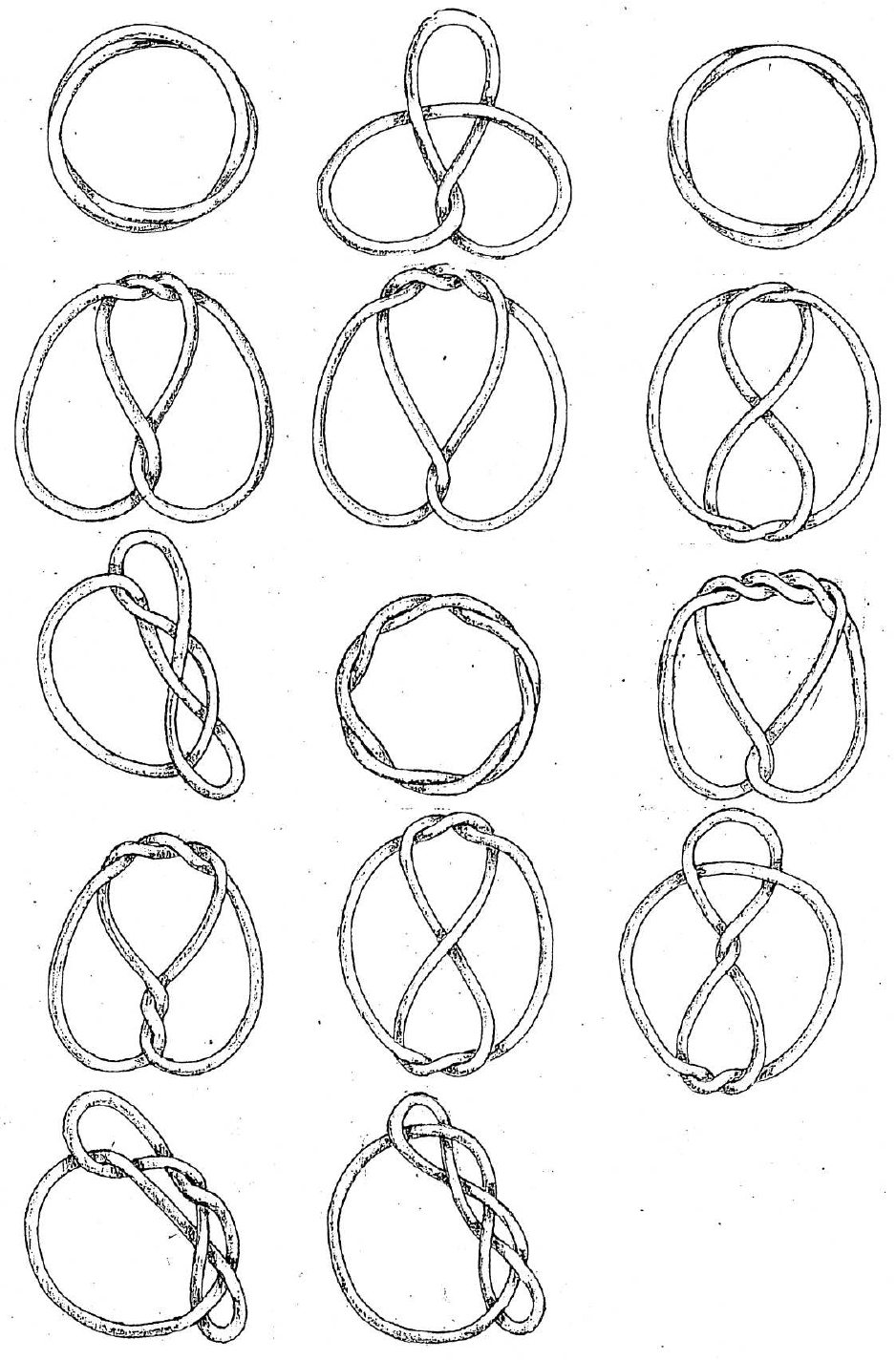}
\par\medskip
{\sc Figure 1.} Mechanical normal forms of prime knots with $7$ crossings or less
\end{figure}

\smallskip 
{\bf Observation 1.} {\sl If our mechanical knot initially forms a knot diagram with seven crossings or less of a prime knot,  then it automatically\pagebreak{} acquires a unique $($up to isometry$)$ mechanical normal form.  These normal forms, which are almost planar, are shown in Figure $1$.}

\smallskip 
All knots with crossing number less than 8 are alternating. Note that amphichieral knots (i.e., knots isotopic to their mirror image, e.g. the $4_1$ knot) actually have {\it two} mechanical normal forms, but they are isometric; we show only one of these forms. Also,  as is customary for knot tables, only one of each pair of nonisotopic mirror symmetric alternating knots (chieral knots) is shown in our table (e.g. only one trefoil, $3_1$).

%

The reader should not misunderstand the meaning of Observation 1: it claims that whenever we {\it started} from a knot diagram {\it with seven crossings or less} 
of a prime knot, our mechanical knot acquired the shape of exactly one of the forms shown in Fig.1. It does not claim, however, that if we start with a knot diagram with, say, seventeen crossings, of one of the prime knots appearing in the knot table with seven crossings or less, then its mechanical normal form will be one of those shown in Figure 1. This is because, first of all, the parameters of our mechanical model (ratio $d/L$ of diameter to length) do not allow us, in practice, to position the device into any knot diagram with seventeen crossings (it is hard enough to do it with seven and sometimes even with six crossings) and, second, because, it is apparently not true (for an idealization of our model with much smaller ratio  $d/L$) that the normal form is unique even for trivial knot diagrams with a very large number of crossings. See the discussion in Sec.4 and Conjecture 2 below. 

Let us note at once that, for a small number of crossings, our mechanical knot outperforms the computer simulation models of knots supplied with an energy functional (see the papers [6], [7]. [8], [9] and the book [1]). Thus, the most popular of the knot energies, M\"obius energy (invented by Jun O'Hara  [9] and studied by Michael Freedman with co-authors [1], [10]) gives two completely different normal forms for the eight knot $4_1$ (shown in Figure 2), whereas our little device, no matter what initial shape of the eight knot we use (including the two shown in Figure 2; note that these knots are not flat at all), immediately takes the shape of the normal form $4_1$ shown in Figure 1; if positioned as the second of the diagrams in Fig.2, it smoothly moves to that normal form, while starting from the first diagram, it first jumps up in the air and acquires that normal form before flopping back on the table. 

\begin{figure}[h]\centering
\includegraphics{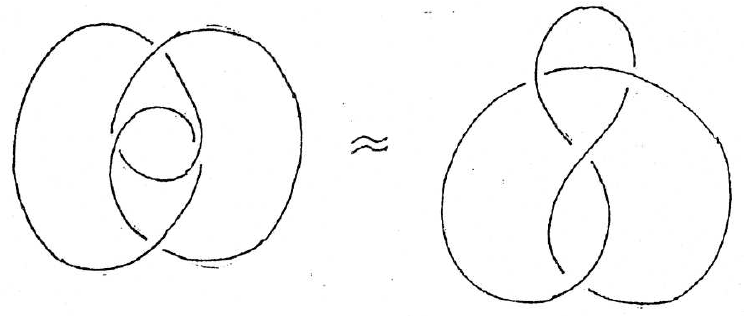}
\par\medskip
{\sc Figure 2.} M\"obius energy normal forms of the eight knot
\end{figure}

\medskip
The reasons for which our mechanical knot sometimes jumps and sometimes only slides are related to Whitney's winding number of plane curves and are discussed below 
(see  Sec.3).

\medskip 
{\bf 2. The idealized models}  

\smallskip
We mentioned ideal versions of our mechanical model above in connection with the question of uniqueness. Here we note that in ideal versions of the model, it should be required that there be no friction preventing the sliding of one branch of the fat knot along another branch, and that the weight of the mechanism be negligible as compared with its resilience (internal energy). In practice, this is not always so: in some cases the surface of the fat knot is not smooth enough for frictionless sliding, and it sometimes requires a little ``prompting" by the experimenter to achieve the normal form; further, for longer mechanical knots (we tested one with $L=70$ cm, $d=2$ mm), the weight of the mechanical knot (with few crossings in its initial position) was sometimes too much for the resilient force of the gadget, which failed to jump off the table and acquire the expected  normal form (unlike a similarly knotted shorter fat knot).

Note that in the idealized versions the ratio $s:=d/L$ must be very small but {\it positive}, and the normal forms depend (up to homothety) on that ratio, so that for idealized mechanical knots we should speak of {\it mechanical $s$-normal forms}. This leads us to the following definition. Let $K$ be a smooth knot of length $L$ in $\R^3$, let $T$ be its tubular neighborhood of cross-section radius $d/2$; position a mechanical knot so as to make it isometric to the tube $T$ and release it; suppose it acquires the normal form 
$T_0$, and denote by $K_0$ the core (central curve) of $T_0$: then $K_0$ is called the 
{\it $s$-normal form of the knot $K$}. 
 
Of course this is not a mathematical definition, since the key ingredient is an unformalized  physical process, but it differs to the previously defined normal form in that it refers to classical knots (rather than fat knots) and is itself a classical knot. To transform this definition into a rigorous mathematical one, we should replace the mechanical device   by an electronic one -- a smooth knot supplied with an appropriate energy functional (depending on $s$), and define the $s$-normal form(s) of the knot as its positions minimizing this functional. Such energy functionals are discussed in Sec.7 
below.  

\medskip 
{\bf 3. The Whitney index and jumps}  

\smallskip
Recall that the {\it winding number} (or {\it Whitney index}) $w(C)$ of a smooth regular curve $C$ in the plane is the number of revolutions that its tangent vector performs when we go around the curve. To every initial position $M$ of our mechanical knot on the table, we can associate a smooth planar curve $C_M$, the projection of the smooth curve in space constituting the core (central curve) of $M$ onto the plane of the table. Without loss of generality, this projection can be assumed to be a smooth immersion (it differs from the corresponding knot diagram in that the crossing points with under/overpasses are replaced by transversal self-intersections). 

Then the winding number $w(C_M)$ will be called the {\it winding number of the mechanical knot $M$}. Numerous experiments with our device led to the following conclusion.

\medskip 
{\bf Observation 2.}  
{\sl As a rule, when the device jumps off the table, a change in the winding number occurs. 
If one free loop disappears (Reidemeister $\Omega_1$), causing the device to jump off the table to acquire its 
normal form, the winding number changes by +1 or -1depending on the orientation of the disappearing loop. 
However, a change in the winding number is not a necessary condition for jumping. For example, the trivial 
knot with two crossing points that consists of two oppositely oriented little loops will jump up of the table 
although its winding number will remain equal to +1.}

\bigskip 
{\bf  4. Local moves of mechanical knots}  

\smallskip
In observing the evolution of the shape of mechanical knots, we were able to notice that the device actually performs certain local moves well known in classical knot theory. In particular, mechanical knots can (under appropriate ``boundary conditions"):   

\smallskip $\bullet$
{\sl perform the first Reidemeister move $\Omega_1$ in both directions}; we have already seen that our device immediately gets rid of little loops, but it can also create a loop: if we take the round unknot and move together two points of the knot, at some moment the knot ``flips", forming a loop (Figure 3);

\begin{figure}[h]\centering
\includegraphics{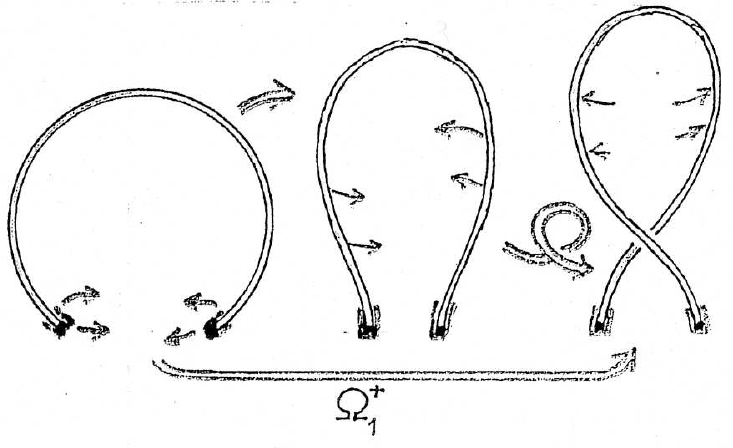}
\par\medskip
{\sc Figure 3.} ``Reidemeister move" forming a little loop
\end{figure}

\smallskip $\bullet$
{\sl perform the second Reidemeister move $\Omega_2$ in both directions} (Figure 4); 
note that the boundary conditions for the simplifying move (in which two crossing points disappear) are not the same as those for the move $\Omega_2^+$ which creates two new crossing points; 
in the second case the ``boundary points" $A$ and $B$ are fixed, in the first  
one they are allowed to slide, but along fixed directions; we also note that, 
under certain conditions, the choice of the branch of the knot that passes 
over the other in the move $\Omega^+_2$ is not predetermined (in a sense, it 
must be considered random), although under other conditions the global configuration of the knot determines which of the branches overpasses;  

\begin{figure}[h]\centering
\includegraphics{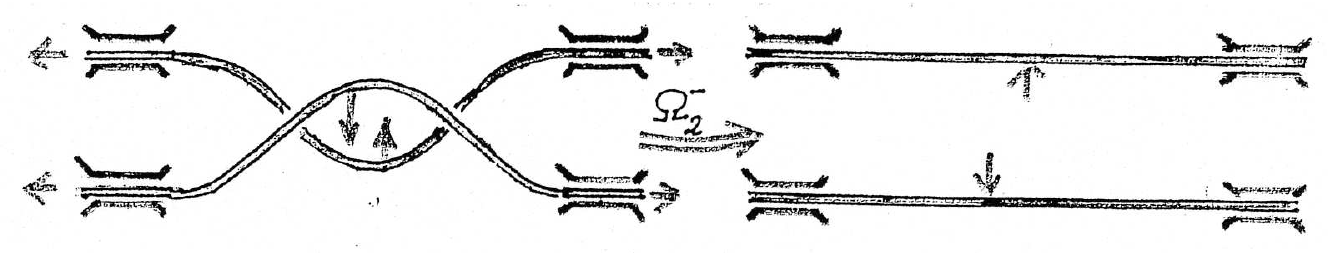}
\par\medskip
{\sc Figure 4.} The second Reidemeister move destroying crossing points
\end{figure}

\medskip $\bullet$
{\sl perform the second Reidemeister move $\Omega_2$} (Figure 5); here the ``boundary points" must be allowed to slide in the directions indicated. 

\begin{figure}[h]\centering
\includegraphics{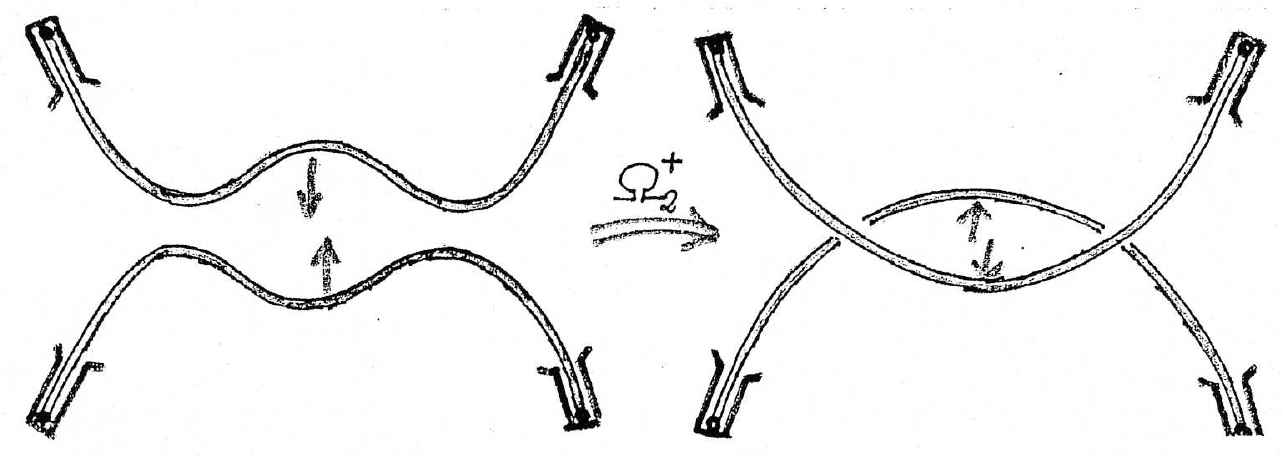}
\par\medskip
{\sc Figure 5.} The second Reidemeister move creating crossing points
\end{figure}

\medskip $\bullet$
{\sl perform the third Reidemeister move $\Omega_3$} (Figure 6); here the ``boundary points" must be allowed to slide in the directions indicated. 
\begin{figure}[h]\centering
\includegraphics{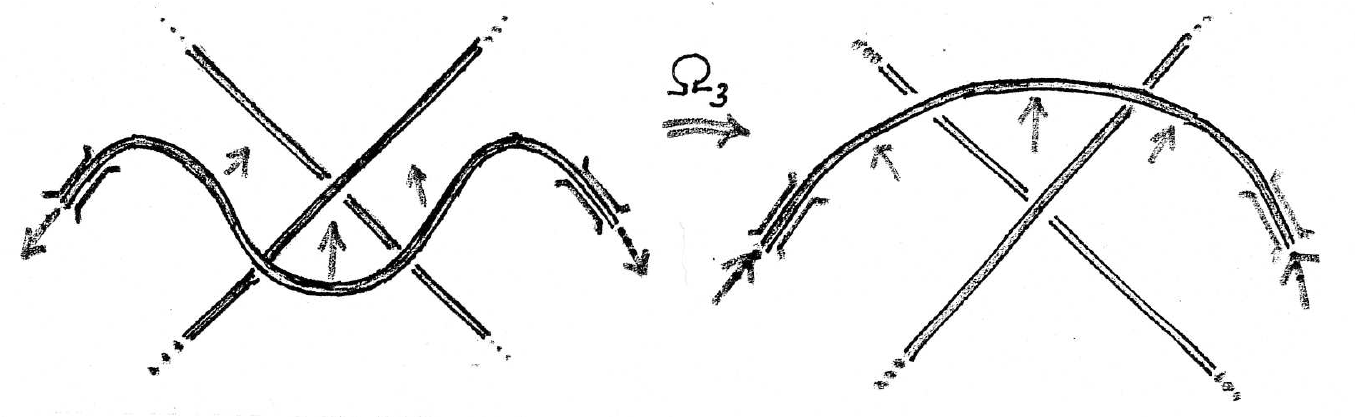}
\par\medskip
{\sc Figure 6.} The third Reidemeister move
\end{figure}

\medskip
Note that, although our mechanical knot can do the moves $\Omega_1$ and $\Omega_2$ in both directions, it ``prefers" performing them in the simplifying directions (decreasing the number of crossing points), because the appropriate\pagebreak{} boundary conditions for these moves occur ``much more often" than those creating new crossing points. 


\medskip
It should also be noted that the Reidemeister moves as performed by our mechanical knot look like the corresponding local moves, but actually are not, rigorously speaking, local. Indeed, when the moves, say,  $\Omega_2^{\pm}$ are performed, the device slides back and forth past the boundary points, thus modifying the global shape of the knot. What is important, however, is not how strictly local the moves are, but what they actually do, namely, in the case of $\Omega_2^{\pm}$, they remove an overlap or add an overlap.  

The mechanical knot can also do the Markov moves, but we have observed such moves only in one direction. Recall that the Markov moves are performed on knots (or links) presented as the closure of a braid; denoting the braid group on $n$ strands by $B_n$ and by $b_1,\dots ,b_{n-1}$ its standard Artin generators, the first move, called {\it stabilization}, has the form  
$$
b \longleftrightarrow bb_{n}^{\pm 1}, \quad \text {where} \quad b\in B_n  
$$
(note that here $b_n$ is a generator of $B_{n+1}$ and does not belong to $B_n$), while the second move is just {\it conjugation},
$$
b \longleftrightarrow a^{-1}ba, \quad \text {where} \quad a,b\in B_n.  
$$

To perform Markov moves, we used an ``almost planar" version of our device that we call  
a flat knot. More precisely, by a {\it flat mechanical knot} we mean the set of different positions of our device placed on the table but constrained from above by a glass pane secured just above the table top and parallel to it; for a mechanical knot of diameter $d=2$mm, the distance between the glass pane and the table was just over 6mm, thus allowing a branch of the knot to move over a crossing (as in an $\Omega_3$ move), but not allowing it to jump up and destroy little loops. This device was able to 

\smallskip $\bullet$
{\sl perform the first Markov move from right to left, i.e., 
$bb_{n}^{\pm 1}\longrightarrow b$}; actually this is just another instance of removing a little loop (Figure 7).

\begin{figure}[h]\centering
\includegraphics{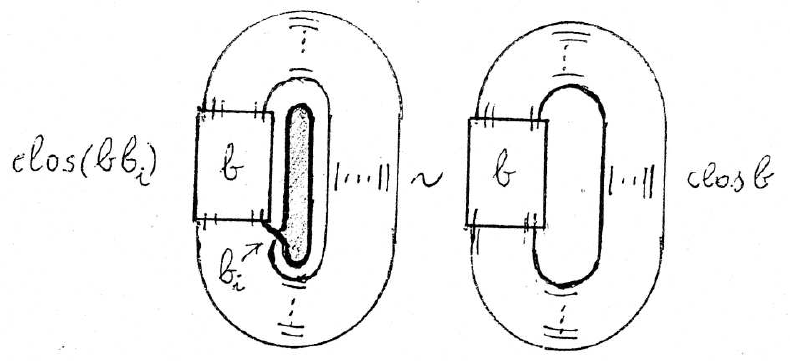}
\par\smallskip
{\sc Figure 7.} First Markov move (destabilization)
\end{figure}

\medskip $\bullet$
{\sl perform the second Markov move from right to left, i.e., 
$a^{-1}ba \longrightarrow b$}; this move is quite spectacular: the braids $a^{-1}$, 
$a$, $b$ begin stretching along the parallel strands of the closure until  $a^{-1}$
meets $a$ and  they cancel each other (as in the $\Omega_2^-$ move), while $b$ continues stretching until it is uniformly distributed along the knot (Figure 8). 

\begin{figure}[h]\centering
\includegraphics{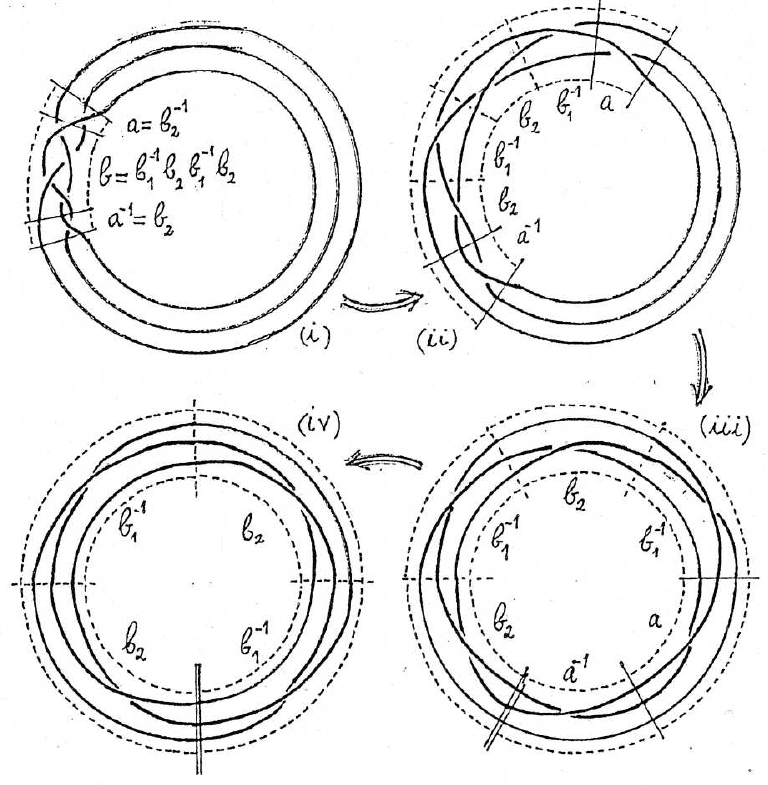}
\par\smallskip
{\sc Figure 8.} Second Markov move (deconjugation)
\end{figure}

Note that in the example shown in the figure, the evolution of the knot does not end there, and in a second phase of the evolution, when the constraining glass pane is removed, the knot eventually acquires the normal form $4_1$ shown in Figure 1. 

\medskip
One of the really remarkable properties of our clever little device is that it can perform the 
{\it Whitney trick} (see e.g. [6], p.17), i.e., mutually destroy two successive little loops  as shown in Figure 8 simply by locally minimizing its curvature.  

\begin{figure}[h]\centering
\vskip-3pt
\includegraphics{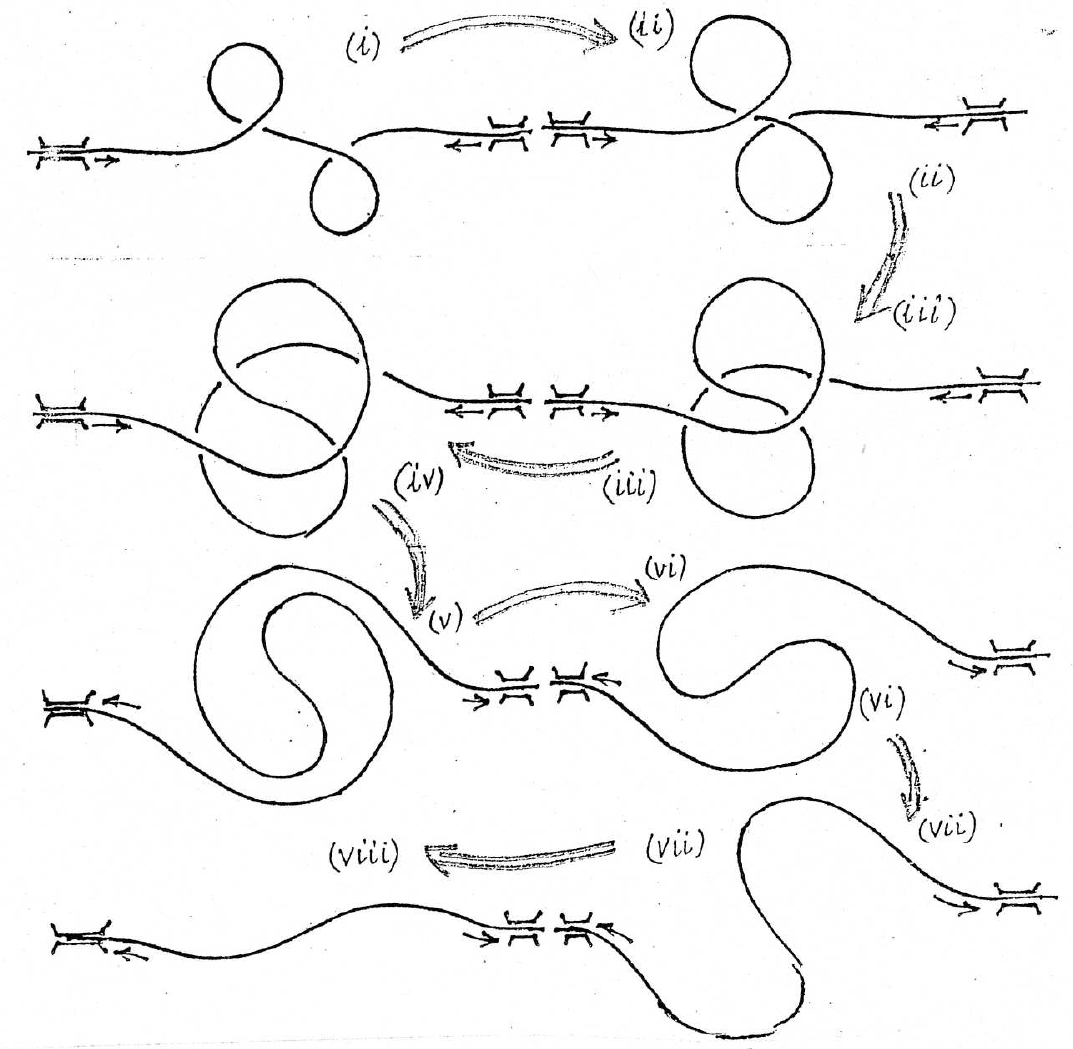}
\par\smallskip
{\sc Figure 8.} The Whitney trick
\vskip-6pt
\end{figure}

\medskip
It should be understood that, in the process depicted in the figure, here again the device was constrained from above by a pane of glass placed just above the table and parallel to it (otherwise the loops would have forced the knot to jump up off the table and each would have self-annihilated itself  independently of the other one).

\smallskip 
{\bf Remark.} The simplest nontrivial knot, the trefoil, can be presented by a smooth curve of total curvature $4\pi$. A beautiful theorem due to John Milnor [11] says that a smooth curve in $\R^3$ with total curvature less than $4\pi$ is unknotted. There is no hope, however, that total curvature could be a measure of the ``complexity" of knots (like their crossing number), because for any positive $\varepsilon$, however small, one can easily construct a knot as complicated as we like (e.g. having a crossing number greater than any apriori given number) of total  curvature less than $4\pi+\varepsilon$. Nevertheless, one might conjecture that total curvature together with some other differential-geometric or topological characteristics of smooth curves might provide a measure of complexity. In this connection, see the discussion in \S 7. 

\bigskip 
{\bf  5. Nonuniqueness of normal forms}  

\smallskip
Experimenting with our mechanical device (of diameter to length ratio $s=d/L=0.004$), we were unable to find any knots that have more than one mechanical normal form, and we claim that no such knots (with 7 crossings or less) exist. However, we believe that uniqueness of mechanical normal forms holds only for a small number of crossing points (less than 30, less than 20?). More specifically, having in mind idealized mechanical models with tiny values of $s=d/L$, we conjecture the following. 

\smallskip 
{\bf Conjecture 1.} {\sl For any $n_0>0$, there exists a sufficiently small $s>0$ such that the unknot has an $s$-mechanical normal form with $n>n_0$ crossing points.}

\smallskip
Recall that for M\"obius energy, the unknot has only one normal form (energy minimum), namely the round circle [3]. 

\smallskip 
{\bf Conjecture 2.} {\sl For any sufficiently large $n_0>0$, there exists a positive $s$ and a nontrivial knot with crossing number $n$ greater than $n_0$ which has more than one $s$-normal form.}

\smallskip
Examples justifying both of these two conjectures can be constructed by the same method: it can be shown that certain different braids have the same closure, but the isotopy equivalence of these closures involves Markov moves that increase as well as decrease the length of the braid, whereas our device can only do the simplifying Markov moves.   

In order to construct such examples, and verify their validity, we need to be able to work with highly idealized mechanical models with tiny values of the parameter $s$, and this can only be done by computer simulation after an appropriate energy functional (imitating the resilience of our device) has been found.  

\pagebreak
{\bf \S 6. Flat knots}  

\smallskip
Flat mechanical knots were described in the previous section. Numerous experiments of the following type were performed with them: the device was knotted (i.e., given the shape of some planar knot diagram and pressed against the table top by a pane of glass) and then the glass pane was raised to the height of approximately 7mm above the table. The device would then move very rapidly on the surface of the table, with some of its branches sliding along other ones, almost instantly acquiring what we call its {\it flat mechanical normal form}.

\smallskip 
{\bf Observation 3.} {\sl For a small number of crossings in its original position, the flat mechanical knot acquires a unique flat normal form.}

\smallskip
Some such normal forms are shown in Figure 9. The ``small number" here is at least 5. There are reasons to believe that it is at least ten or even more, but experiments were carried out systematically only for five crossings or less, and the total number of experiments was several hundred (rather than many thousands for spacial mechanical knots). 

\begin{figure}[h]\centering
\includegraphics{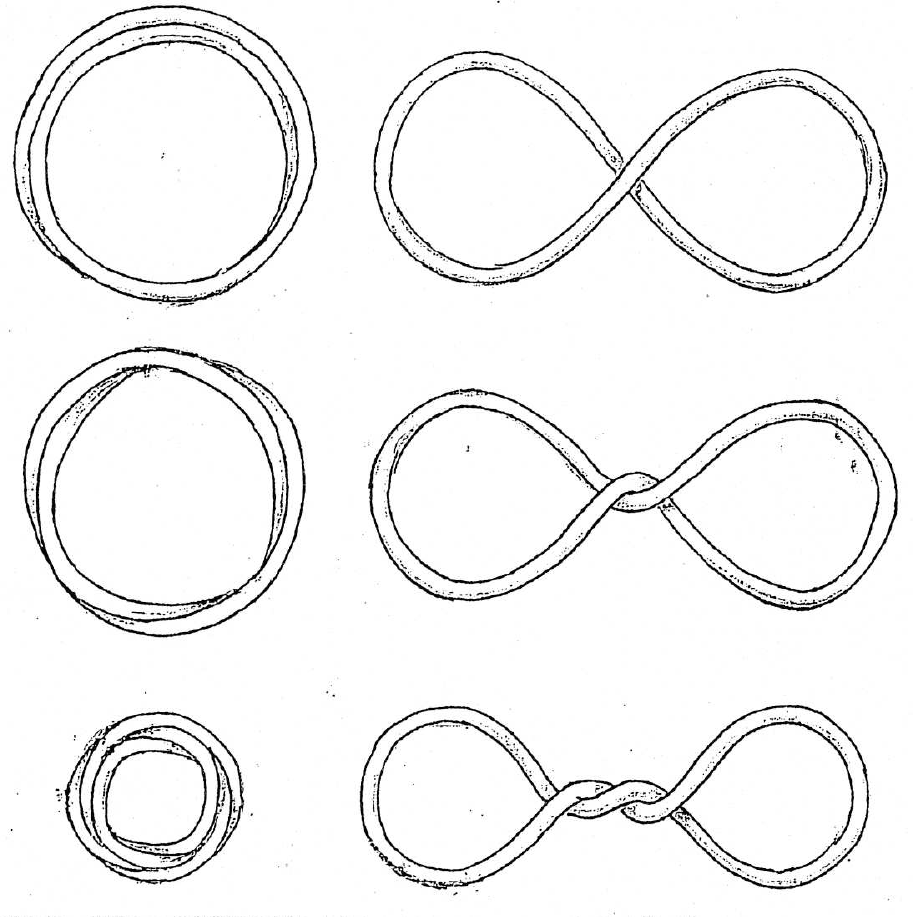}
\par\medskip
{\sc Figure 9.} Some flat normal forms
\end{figure}

\medskip
Note that if  the initial position of the flat knot has little loops, then these loops sometimes become incorporated into the normal form (increasing in size, as for the second normal form in Figure 8), sometimes mutually cancel each other as in the Whitney trick (see Figure 7).

The physical notion of flat mechanical knot has a mathematical counterpart: a 
{\it flat $($mathematical\,$)$ knot} can be defined as the equivalence class of all knot diagrams under the Reidemeister moves $\Omega_2$ and $\Omega_3$ (the destruction and creation of little loops, $\Omega_1$, is not allowed). When two knot diagrams are equivalent in this sense, we will say that they are {\it flat isotopic}. Flat isotopy is a stronger equivalence relation than the usual ambient isotopy relation, nevertheless we have the following theorem. 

\medskip {\bf Theorem 1.}
{\sl If two ambient isotopic knot diagrams $K_1$ and $K_2$ have the same Whitney index, then they are flat isotopic. }

\medskip {\it Proof  .} 
By the Reidemeister lemma, there exists a sequence of Reidemeister moves  taking $K_1$ to $K_2$. We will replace this sequence by another sequence, not involving any $\Omega_1$ moves, that nevertheless takes $K_1$ to $K_2$, in the following way. 

If an $\Omega_1$ move that destroys a little loop is performed at some point, we do not perform it, but, instead, tighten the little loop until it becomes tiny (we call it a {\it kink}) and ``freeze" it to the part of the diagram where it is located. During further moves of the diagram, the kink stays frozen to the diagram and moves rigidly with it as if it were simply a marked point on the diagram. (Actually, when that part of the diagram participates in  moves where it crosses other parts of the diagram, the kink moves along with it without changing its shape,  but as it does this, tiny $\Omega_2$ and $\Omega_3$ moves involving it occur.) 

If an $\Omega_1$ move that creates a little loop is performed at some point, we do not perform this move, but carry out the Whitney trick (which only uses $\Omega_2$ and $\Omega_3$ moves) to create two opposite little loops; one of them plays the role of the loop created by the $\Omega_1$ move, while the other is tightened and becomes a frozen kink (and we forget about it until the end of the construction). 

Proceeding in this way, after a finite number of $\Omega_2$ and $\Omega_3$ moves 
(the old ones from the original sequence and the new ones -- those used to mimic the 
$\Omega_1$ moves that create new loops and those appearing when the kinks cross other parts of the knot diagram), we will obtain a knot diagram $K_2'$ which is exactly like $K_3$ except that it has a certain number of kinks (little loops) on it. The knot diagram $K_2'$ has the same Whitney index as $K_1$ (because $\Omega_2$ and $\Omega_3$ moves do not change the Whitney index) and so, by assumption, has the same Whitney index as $K_2$. Thus the little loops on $K_2'$ destroy each other and we obtain $K_2$ without doing any $\Omega_1$ moves, i.e., by means of a flat isotopy. 
$\hphantom{A}$ \hfill{$\square$}

\medskip
Flat knots $K$ have two well-known invariants, the Kauffman bracket $\langle K \rangle$ (see, e.g., [5, p.23-24]) and the Whitney index or winding number $w(K)$ (see Sec.3  above). 

To each flat knot $K$, one can assign another integer invariant, its {\it twisting number} $\tau(K)$, defined as shown in Figure 10: given a knot diagram $K$, we first transform it into an unknot by appropriate crossing changes (see, e.g., [5], Chapter 1), then replace that unknot by a flat ribbon, cut the ribbon and pull it taut (this is the well known belt trick, see, e.g., [5]); then the number of twists obtained (with sign specified via the ``left-hand rule") is, by definition, $\tau(K)$. As far as I know, this definition first appeared in S.Matveev's article [12] (see also [13]) for the particular case of the flat unknot. 

\begin{figure}[h]\centering
\includegraphics{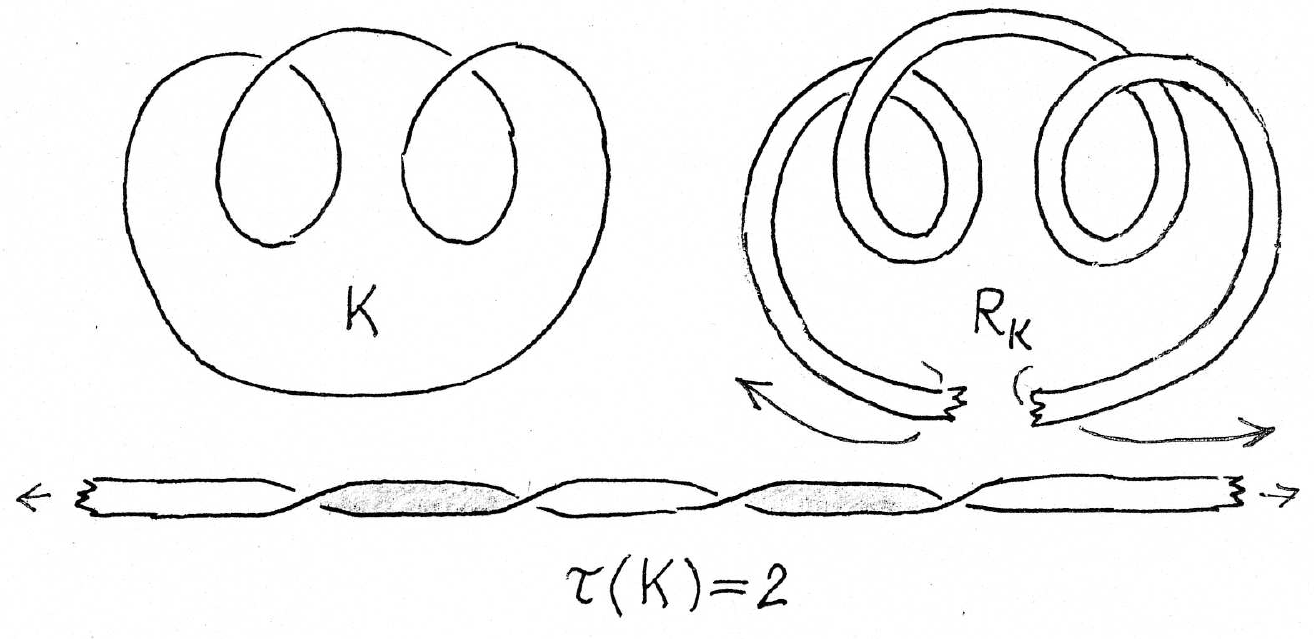}
\par\medskip
{\sc Figure 10.} Definition of the twisting number
\end{figure}

\medskip 
{\bf Theorem 2.} {\sl Flat knots possess the following  three invariants: the Kauffman bracket $\langle K \rangle$, the Whitney index $w(K)$, and the twisting number $\tau(K)$. Together, $w(K)$ and $\tau(K)$ constitute a complete system of invariants for flat unknots.}   

\smallskip {\it Proof.}
The invariance of the Kauffman bracket is a classical lemma due to Kauffman 
(see, e.g., [5], pp.24-25). To prove the invariance of the twisting number, it suffices to check that it does not change under $\Omega_2$ and $\Omega_3$, but this is a straightforward verification. The invariance of the Whitney index, as well as the fact that together with the twisting number it constitutes a complete system of invariants for flat unknots, was proved by Matveev in [12] by an elementary but tedious argument (see also [13]). 
$\hphantom{A}$ \hfill{$\square$}

\medskip 
Of course the algebraic classification problem for all flat knots is more difficult than the same problem for classical (spacial) knots and would require much stronger invariants (of Vassiliev type?).

\bigskip 
{\bf \S 7. Perspectives: new knot energies for computer simulations}  

\medskip
The analytical specification of a knot energy mimicking the resilience of our little 
device is ongoing joint research with Oleg Karpenkov and will be published elsewhere.

Here we only mention that our functional is defined for flat knots only, it 
depends on a small parameter $d$ in the space of smooth flat knots in $\R^2$; roughly speaking, the functional is based on an ``unbending force" (the larger the curvature in the vicinity of some point, the stronger the force that tends to straighten out the curve near that point). For flat knots, unlike ``fat knots" of cross section diameter $d$ in space, no  $\delta$-function type of repulsive term that does not allow different branches of the curve to be at a distance {\it strictly smaller} than $d$ from each other, but allows theses branches to be at a distance {\it exactly equal to $d$}, is needed. 

However, because of Theorem 1, an energy functional for flat knots suffices to determine if two knot diagrams correspond to the same knot (i.e., are ambient isotopic in space). To do this, {\sl compute the Whitney index of the two knot diagrams, add an appropriate number of little loops to one of them so as to equalize their Whitney indices and then apply the energy functional; the given knot diagrams present the same knot if and only if    
the two flat normal forms obtained are identical.}

Of course this works only provided that the number of crossings is small enough to ensure uniqueness of normal forms. An optimistic prediction is that ``small" means 
less than 30; if this were so, we would have a fast effective algorithm for comparing knot diagrams with 30 crossings or less. 

To my mind, a more interesting problem is to define other energy functionals that do not necessarily mimic our little device. Such functionals should include not only a $\delta$-function-like repulsive term, but other terms related to curvature, torsion, twisting numbers, and whatnot.    

More specifically, such a functional could mimic another concrete device with which I have performed a series of similar experiments. This new device is made from a 
thin ($\sim 3$mm) three-strand metallic cable which resists twisting (unlike the device described in this article, it does the ``belt trick"). The corresponding functional would involve three summands: an unbending term (related to local curvature), a repulsive  
term (related to the thickness of the device), and an untwisting term (related to torsion and/or the so called ``blackboard framing" of the initial knot diagram).  

\bigskip
\centerline{{\sc Acknowledgements}}

\medskip 
I am is grateful to Elena Efimova, for providing me with the device used in the experiments (and for many other things) and to Mikhail Panov and Victor Shuvalov, who helped with the illustrations. I also ackowledge the support of the RFBR-CNRS-a grant 
$\#$10-01-93111.    

\bigskip
\centerline{{\sc References}}

\medskip 
[1] Jun O'Hara, {\it Energy of Knots and Conformal Geometry}, World Scientific, Singapore, 2003

\smallskip
[2] M.H.Freedman, Z.-X.He, Z.Wang, M\"obius energy of knots and unknots,  Ann. of  Math. (2) {\bf 139}, 1-50, (1994)

\smallskip
[3] K.Ahara, Energy of a knot. video available on {\it youtube.com} 

\smallskip
[4] Ying-Qing Wu, Polygonal knot energy, available at 
    
\centerline{ {\it  http://www.uiowa.edu/~wu/min/ming.html }}  

\smallskip
[5] V.V.Prasolov, A.B.Sossinsky,  Knots, Links, Braids and 3-Manifolds, AMS Publ., 
Providence, R.I. 1997

\smallskip
[6] H.K.Moffat {\it The degree of knottedness of tangled vortex lines} J.Fluid Mech.
{\bf 35} (1), (1969) 117-129

\smallskip
[7] V.I.Arnold {\it Asymptotic Hopf invariant and its applications} [In Russian] In the book: Proceedings of the All-Union School in Differential Equations, Erevan, Acad. Sc. Armenian SSR (1974) 229-256. 

\smallskip
[8] H.K.Moffat {\it Magnetostatic equilibria and analogous Euler flows of arbitrarily complex topology}, J.Fluid Mech. {\bf 159}, (1985) 359-378.

\smallskip
[9] Jun O'Hara, Energy of a knot, Topology {\bf 30}, 241-247 (1992)  

\smallskip
[10] M.H.Freedman, Z.-X.He, Links of tori and the energy of incompressible flows, Topology, {\bf 30}, 283-387, (1992)

\smallskip
[11] J.W.Milnor {\it On the total curvature of knots}, Ann. of Math. {\bf 52} (1950), 
248-257.

\smallskip
[12] S.V.Matveev, Unraveling closed curves in the plane  [in Russian], Kvant, 1983, 
no.4, 22-28. 

\smallskip
[13] A.Sossinsky, Wire loops on the plane [in Russian], in: Problems of the Tournament of Towns, MCCME Publishers, Moscow, 2009

\end{document}